\documentclass[letterpaper, 10 pt, conference]{ieeeconf}  

\IEEEoverridecommandlockouts                              

\usepackage{algorithm2e}
\usepackage{algorithmicx}

\usepackage{tikz}
 \usetikzlibrary{arrows,shapes,positioning}
 \usetikzlibrary{decorations.markings}
 \tikzstyle arrowstyle=[scale=1]
 \tikzstyle directed=[postaction={decorate,decoration={markings,mark=at position .65 with {\arrow[arrowstyle]{stealth}}}}]
 \tikzstyle reverse directed=[postaction={decorate,decoration={markings,mark=at position .65 with {\arrowreversed[arrowstyle]{stealth};}}}]
\usepackage{amsmath}
\usepackage{epstopdf}
\usepackage{latexsym, amssymb,booktabs}
\usepackage{graphicx}
\usepackage{setspace}
\usepackage{amsmath, amsbsy}
\usepackage{amsopn, amstext}
\usepackage{hyperref,url}
\usepackage{floatrow}
\newfloatcommand{capbtabbox}{table}[][\FBwidth]
\usetikzlibrary{arrows}
\DeclareMathOperator{\rank}{rank}

\DeclareMathOperator{\Span}{Span}
\DeclareMathOperator{\Col}{Col}
\DeclareMathOperator{\Row}{Row}

\DeclareGraphicsExtensions{.eps,.pdf}

\def\ra{\rightarrow}

\def\b{\beta}
\def\d{\delta}

\def\D{\Delta}

\newcommand{\R}{{\mathbb R}}

\newtheorem{thm}{Theorem}[section]
\newtheorem{lem}[thm]{Lemma}
\newtheorem{cor}[thm]{Corollary}
\newtheorem{dfn}[thm]{Definition}

\newtheorem{exa}[thm]{Example}
\newtheorem{rem}[thm]{Remark}
\newtheorem{alg}[thm]{Algorithm}

\title{Potential Games Design Using Local Information}
\author{Changxi Li$^{1}$,  Fenghua He$^{1}$, Hongsheng Qi$^{2}$, and Daizhan Cheng$^{2}$
\thanks{*This work is supported partly by the National Natural Science Foundation of China (NSFC) under Grants 61473099, 61773371, 61733018 and 61333001.}
\thanks{$^{1}$ Changxi Li and Fenghua He  are with Harbin Institute of Technology, Harbin 150001, P.~R.~China {\tt\small  changxi1989@163.com, hefenghua@hit.edu.cn}}
\thanks{$^{2}$ Hongsheng Qi and Daizhan Cheng are with Key Laboratory of Systems and Control,
		Academy of Mathematics and Systems Sciences, Chinese Academy of Sciences,
		Beijing 100190, P. R. China {\tt\small  qihongsh@amss.ac.cn, dcheng@iss.ac.cn}}
\thanks{ Corresponding author: Fenghua He. Tel.: +86 0451-86402947; fax.: +86 0451-86414580.}
}

\begin{document}

	\maketitle
	\thispagestyle{empty}
	\pagestyle{empty}


\begin{abstract}
Consider a multiplayer game, and assume a system level objective function, which the system wants to  optimize, is given. This paper aims at accomplishing this goal via potential game theory when players can only get part of other players' information. The technique  is designing a set of local information based utility functions, which guarantee that the designed game is potential, with the system level objective function its potential function.
 First,   the existence of local information based utility functions can be verified by checking whether the corresponding  linear equations have a solution.
Then an algorithm is proposed to calculate the local information based utility functions when the utility design equations have  solutions.
Finally,  consensus problem of multiagent system is considered to demonstrate the effectiveness of the proposed design procedure.
\end{abstract}


\section{Introduction}
Game-theoretical control  has drawn considerable attention in recent years due to its widespread applications. Some representative works include: (i) consensus/synchronization of multi-agent systems \cite{jm09}; (ii) distributed optimization \cite{li13, by10}; (iii) control in wireless networks; (iv) optimization in energy \cite{wss12} and transportation networks \cite{pn17, wan13}, just to name a few.

The  content of game-theoretical control is using game theory to solve control problems in interacting setting, such as multiagent systems \cite{li13, by10}.
Addressing such issues via game theory needs two steps. The first step is to view the agent as an intelligent rational decision-maker in a game with defining a set of  available actions and utility function for every player. The second step is to specify a learning rule for the designed game so that the agents can reach a desirable situation, e.g., a Nash equilibrium. Therefore, there are two basic  tasks in game-theoretic control: utility design and learning rule design \cite{RG11}. Utility functions describe  components' incentives, and learning rules mean how  each player processes its available information to formulate a decision.

\begin{figure}[!hbt]
    \centering
    \includegraphics[height = 2.3cm, width =3.2 cm]{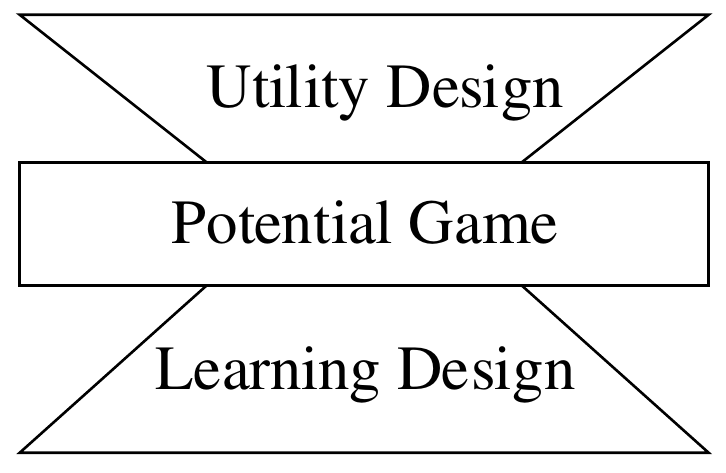}
    \caption{Hourglass Architecture of Game-Theoretical Control}
    \label{Fig1}  
\end{figure}

Compared with traditional  methods, the advantage of game-theoretical control is that it provides a modularized design architecture, i.e. we can design utility functions and learning rules separately \cite{jr15}. The separation is described as an hourglass architecture in \cite{RG11}, which is shown in  Fig. \ref{Fig1}.
When designing games, one idea is to make sure that the designed game falls under some special category games, such as potential games \cite{mon96}.  One advantage of designing the game as potential game is that there are a variety of learning rules which lead to a Nash equilibrium, e.g. myopic best response, log-liner learning, and fictitious play \cite{ca10}-\cite{js05}. Several papers have devoted to potential game based design in distributed control \cite{ec08}-\cite{jra09}. Other game design methods include wonderful life utility design \cite{jra12}, Shapley value utility design \cite{ea08}, congestion game based design \cite{hy17}, etc. However most of the above works provide no systematic methods on designing local information based utility functions. Here  local information means that players can only get part of other players' information when they play the designed game, such as  networked game.

This paper focuses on providing a systematic method for designing finite potential game using local information. As far as we know, the most relevant works are \cite{li13} and \cite{liu16}. But our work is different to theirs. \cite{li13} provided a systematic methodology for designing potential games with continuous action sets, where the utility functions are local information based. It showed that for any given system level objective function, there exists at least one method to design the local information based utility functions \cite{li13}. However, when we turn to games with finite  action sets, the existence is not guaranteed. As for \cite{liu16}, it  presented a necessary and sufficient condition for the existence of local information based utility functions. But no systematic method is provided for designing local information based utility functions. Furthermore, the learning rule used in  \cite{liu16} is better reply. Using better reply, local information based potential game can  converge to a Nash  equilibrium, but may not  a maximum point of the system objective function.

The contributions of this paper are threefold:
(i) A necessary and sufficient condition for the existence of local information based utility functions is obtained, which can be verified by  checking whether the corresponding linear equations have  solutions.
(ii) A  method for designing finite potential game using local information is presented when the linear equations have a solution.
(iii) An example on consensus problem is provided to demonstrate the effectiveness of the design procedure.

The rest of this paper is organized as follows: Section II provides some preliminaries, including semi-tensor product (STP) of matrices, game theory, and problem description. Section III considers the design of potential game using local-based information. Section IV considers  application of the design method to consensus problem. A brief conclusion is given in Section V.

\textbf{Notations}: $\R^n$ is denoted by the Euclidean space of all real $n$-vectors. ${\cal M}_{m\times n}$ is the set of $m\times n$ real matrices.
${\bf 1}_{\ell}=(\underbrace{1,1,\cdots,1}_{\ell})^T$. $I_m\in {\cal M}_{m\times m}$ is the $m\times m$-dimensional identity matrix.
${\bf 0}_{m\times n}\in {\cal M}_{m\times m}$ is the $m\times m$-dimensional zero matrix.
${\cal D}_k:=\left\{1,2,\cdots,k\right\}, k\geq 2$. $\d_n^i$ is denoted by the $i$-th column of the identity matrix $I_n$. $\Col(M)$ ($\Row(M)$) is the set of columns (rows)  of $M$. The transposition of matrix $A\in{\cal M}_{m\times n}$ is denoted by $A^\mathrm{T}\in {\cal M}_{n\times m}.$  $\Span\{V_1,\cdots,V_s\}$ is the subspace spanned by $\left\{V_i\;|\;V_i\in\R^n,i=1,\cdots,s\right\}$.

\section{Preliminaries}
\subsection{Semi-tensor Product of Matrices}
The basic tool used in this paper is STP of matrices. We give a brief survey on STP of matrices. Please refer to \cite{che12} for  more details.

\begin{dfn} \label{da.1} \cite{che12}
Suppose  $A\in {\cal M}_{m\times n}$, $B\in {\cal M}_{p\times q}$, and $l$ be the least common multiple of $n$ and $p$.
The STP of $A$ and $B$ is defined by
\begin{align*}
A\ltimes B:= \left(A\otimes I_{l/n}\right)\left(B\otimes I_{l/p}\right)\in {\cal M}_{ml/n\times ql/p},
\end{align*}
where $\otimes$ is the Kronecker product.
\end{dfn}

Assume $i\in\mathcal{D}_k$. By identifying $i\sim\delta_k^i$ we call $\delta_k^i$ the vector form of  integer $i$. A function $f:\prod_{i=1}^n{\cal D}_{k_i}\ra \R$ is called a mix-valued pseudo-logical function.
\begin{dfn}\cite{che12}
Let  $f:\prod_{i=1}^n{\cal D}_{k_i}\ra \R$ be a mix-valued pseudo-logical function. Then there exists a unique row vector $M_f\in \R^{k}$, such that
$$
f(x_1,\cdots,x_n)=M_f\ltimes_{i=1}^n x_i.
$$
$M_f$ is called the structure vector of $f$, and $k=\prod\limits_{i=1}^nk_i$.
\end{dfn}
\subsection{Potential Game}
A finite non-cooperative game is a triple $G=\left\{N,\{S_i\}_{i\in N},\{c_i\}_{i\in N}\right\}$, where $N=\{1,2,\cdots,n\}$ is the set of players, $S_i=\{1,2,\cdots,k_i\}$ is the set of strategies of player $i$ for every $i\in N$, and $c_i:S\ra \R$ is the utility function of player $i$, with $S:=\prod_{i=1}^nS_i$ being the strategy profile of the game. Let $S^{-i}:=\prod_{j\neq i}S_j$ be the set of partial strategy profiles other than player $i$. Denote ${\cal G}_{[n;k_1,\cdots,k_n]}$ by  set of  form finite games with $|N|=n$, $|S_i|=k_i$, $i=1,\cdots,n$.

 Using the vector expression of strategies, the utility function can be expressed as
$$c_i(x)=V_i\ltimes_{j=1}^nx_j,$$
where $x_i\in S_i$, and $V_i\in \R^k$ is called the structure vector of $c_i$,  $k=\prod_{i=1}^nk_i$.

The concept of potential game was firstly proposed by Rosenthal \cite{ros73}, whose definition is as follows:
\begin{dfn} \label{d2.4} A finite game $G\in {\cal G}_{[n;k_1,\cdots,k_n]}$ is a  potential game if there exists a function $P: S\rightarrow \R$,  such that for every player $i\in N$ and every $s^{-i}\in S^{-i}$
\begin{align*}
\begin{array}{cl}
&c_i(x_i,s^{-i})-c_i(y_i, s^{-i})=P(x_i,s^{-i})-P(y_i, s^{-i}),\\
&~~~~~~~~~~~~~~~~~~~~~~~~~~~~~~~~~~~~~~~\forall x_i,y_i\in S_i,
\end{array}
\end{align*}
where $P$ is called the potential function of  $G$.
\end{dfn}

The following Lemma is obvious according to Definition \ref{d2.4}.
\begin{lem}\label{lem2.3}\cite{che14}
A finite game $G\in{\cal G}_{[n;k_1,\cdots,k_n]}$ is   potential  if and only if there exist functions $d_{i}: S^{-i}\rightarrow \R,i\in N$ such that for every $x\in S$
\begin{align}\label{2.2}
\begin{array}{cl}
P(x)=c_i(x)-d_i( x^{-i}),~\forall i\in N,
\end{array}
\end{align}
where $P(x)$ is the potential function, and $x^{-i}\in S^{-i}$.
\end{lem}
\subsection{Problem Setup}
Consider a multi-player game $G\in{\cal G}_{[n;k_1,\cdots,k_n]}$  played on the network, which is called networked game (NG).   In fact every player can only obtain its neighbours' information when playing game $G$. The neighbors of player $i$ is denoted by $U(i)$, which is defined as follows,
$$U(i)=\{j\in N|\text{$i$ can communicate with $j$}\}.$$
Suppose a system level objective function $\phi(x)$ is given, where $x\in S$. The system wants to optimize the objective function, while all players can only obtain its neighbour’s information. The optimization problem can be described as
\begin{align*}
\begin{array}{cll}
\max&~&~\phi(x_1,x_2,\cdots,x_n)\\
\text{s.t.}&~&c_i(x)=c_i(x_{U(i)},x_i), \\
~&~&~x_i\in S_i, \forall i\in N,
\end{array}
\end{align*}
where $x_{U(i)}=\{x_j\}_{j\in U(i)}.$

To solve this optimization problem,   the idea is to design a potential game  in which the utility function of every player only depends on its neighbours' information, and the potential function of the designed  game is $\phi(x)$. Then using proper  learning algorithm, such as logit learning \cite{ca10},   players' behavior converges to a  strategy profile that maximize the objective function.

\section{Potential Game Design Using Local-based Information}

\subsection{Utility Design Using Local Information}
We consider the design of potential game using local information in this subsection.
Before designing the game, an operator $\Gamma_{U}$, called the $U$-drawing matrix, is necessary, where $U\subset N$ is a group of players in the game $G\in{\cal G}_{[n;k_1,\cdots,k_n]}$.
Set
$$
\Gamma_{U}:=\otimes_{i=1}^{n}\gamma_i
$$
where
$$
\gamma_i:=
\begin{cases}
I_{k_i},~i\in U\\
{\bf 1}_{k_i}^T,~\text{otherwise}.
\end{cases}
$$
The  $U$-drawing matrix is used to ``draw" the strategies of players in $U$ from $N$ \cite{liu16}:
$$
\ltimes_{j\in U}x_j=\Gamma_{U}\ltimes_{j=1}^nx_j,
$$
where $x_i\in S_i$ is the strategy of player $i$. Particularly, $\Gamma_{-i}:=\Gamma_{N\backslash\{i\}},$ where ${N\backslash\{i\}}$ is the set of players except player $i$.

\begin{thm}\label{th3.1}
Consider a utility-adjustable networked game $G\in{\cal G}_{[n;k_1,\cdots,k_n]}$ with objective function $\phi(x)$
$$\phi(x)=V^{\phi}\ltimes_{j=1}^nx_j.$$
Then  local information based utility function can be designed if and only if all the following equations have a solution
\begin{align}\label{3.3}
\begin{array}{ccl}
T_i\cdot\xi_i=(V^{\phi})^\mathrm{T},
\end{array}
\end{align}
where $T_i=[\Gamma_{N_i}^\mathrm{T},\Gamma_{-i}^\mathrm{T}]$, $\xi_i=[(\xi_i^1)^\mathrm{T}, (\xi_i^2)^\mathrm{T}]^\mathrm{T}$, $\xi_i^1\in\R^{k_{N_i}}$, $\xi_i^2\in\R^{k_{-i}}$, $k_{N_i}=\prod_{j\in N_i}k_j$, $k_{-i}=\prod_{j\neq i}k_j$, and $N_i=U(i)\cup \{i\}$, $\forall i \in N$.

Moreover if the solution $\xi_i, \forall i\in N$ exists,   the local information based utility function  of player $i$ is
\begin{align}\label{3.4}
\begin{array}{ccl}
c_i(x)=(\xi_i^1)^\mathrm{T}\Gamma_{N_i}\ltimes_{j=1}^nx_j, \forall i\in N.
\end{array}
\end{align}
\end{thm}
\noindent{\it Proof:} If the utility function of the networked game $G$ is local information-based, then we have
$$c_i(x)=V_i\ltimes_{j\in N_i}x_j=V_i\Gamma_{N_i}\ltimes_{j=1}^nx_j.$$

Rewrite (\ref{2.2}) into vector form, and substituting  (\ref{3.4}) into (\ref{2.2})  yields
\begin{align*}
\begin{array}{cll}
V^\phi\ltimes_{l=1}^nx_l&=&V_i\ltimes_{j\in N_i}x_j-V_{i}^d\ltimes_{j\neq i}x_j,\\
&=&V_i\Gamma_{N_i}\ltimes_{j=1}^nx_j-V_{i}^d\Gamma_{-i}\ltimes_{j=1}^nx_j,
\end{array}
\end{align*}
where $V^\phi$ and $V_{i}^d$ are the structure vectors of $\phi$ and $d_{i}$, respectively.

Since $x_i\in\D_{k_l},i\in N$ are arbitrary,  we have
\begin{align}\label{3.6}
\begin{array}{cl}
V^\phi=V_i\Gamma_{N_i}-V_{i}^d\Gamma_{-i},\forall i\in N,
\end{array}
\end{align}
which is equivalent to (\ref{3.3}).

If the solution $\xi_i$ exists, then we have
$$V_i=(\xi_i^1)^\mathrm{T},~\forall i\in N,$$
which verifies equation (\ref{3.4}).

\hfill $\Box$

Equations (\ref{3.3}) are called utility design equations, and $T_i, \forall i\in N$ are called utility design matrices.
\subsection{Solving Utility Design  Equations}
In this subsection we first explore some properties of utility design equations. Then   we design an algorithm to calculate  solutions of  utility design equations equations when the solutions exist.

Set
\begin{align*}
\begin{array}{cl}
\Lambda_{N_i}:=\mathop\otimes\limits_{j\in N_i}\Lambda_j, ~~\Upsilon_{N_i}:=\mathop\otimes\limits_{j\neq i}\Upsilon_j,
\end{array}
\end{align*}
where
$$
\Lambda_j:=
\begin{cases}
{\bf 1}_{k_j},~j=i\\
I_{k_j},~j\in U(i)
\end{cases},
~
\Upsilon_j:=
\begin{cases}
{\bf 1}_{k_j},~j\notin N_i\\
I_{k_j},~j\in U(i)
\end{cases}.
$$

Denote by $H_i=[\Lambda_{N_i}^\mathrm{T},-\Upsilon_{N_i}^\mathrm{T}]^\mathrm{T}, \forall i\in N.$ Let $\mathcal{H}_i=\Span\Col(H_i), \forall i\in N.$
\begin{thm}\label{th3.2}
Let $T_i$ be the utility design matrix of player $i$ in (\ref{3.3}). Then
\begin{align}
\Span \Row(T_i)=\mathcal{H}_i^\bot, \forall i\in N.
\end{align}
In other words, every column of $H_i$ is a solution of $T_ix=0.$ And
$$\rank(T_i)=k_{N_i}-k_{U(i)},$$
where $k_{N_i}=\prod_{j\in N_i}k_j$, and $k_{U(i)}=\prod_{j\in U(i)}k_j$
\end{thm}
\noindent{\it Proof:}
We omit the detailed proof due to the space limitation.

\hfill $\Box$

Using Theorem \ref{th3.1} and  Theorem \ref{th3.2}, we have the following results:
\begin{cor}
Consider a utility-adjustable networked game $G\in{\cal G}_{[n;k_1,\cdots,k_n]}$ with objective function $\phi(x)$
$$\phi(x)=V^{\phi}\ltimes_{j=1}^nx_j.$$ Then the following four statements are equivalent.
\begin{description}
  \item[(i)]  The local information based utility functions can be designed.
  \item[(ii)] The following equations have a solution
$$T_i\cdot\xi_i=(V^{\phi})^\mathrm{T},\forall i\in N.$$
  \item[(iii)] Suppose $T_i$ is defined in (\ref{3.3}). Then,
  \begin{align}\label{3.7}
  \rank[T_i,(V^{\phi})^\mathrm{T}]=k_{N_i}-k_{U(i)},\forall i\in N.
  \end{align}
  \item[(iv)]
  $$(V^{\phi})^\mathrm{T}\in \mathop\bigcap\limits_{i=1}^n\Span\Col(T_i).$$
\end{description}
\end{cor}
\noindent{\it Proof:}
(i)$\Leftrightarrow$(ii): It is obvious using Theorem \ref{th3.1}.

(ii)$\Leftrightarrow$(iii): From linear algebra we know that condition (\ref{3.7}) is the
necessary and sufficient condition for equation (\ref{3.3}) to have solutions.

(iii)$\Leftrightarrow$(iv): From Theorem \ref{th3.2} one sees that
$$\rank(T_i)=k_{N_i}-k_{U(i)}.$$
Combining equation (\ref{3.7}) we have the following result
$$(V^{\phi})^\mathrm{T}\in \Span\Col(T_i), \forall i\in N,$$
which is equivalent to statement (iv). It is obvious that (iv) implies (iii).

\hfill $\Box$

In the following we design an algorithm to calculate the local information based utility function for each player when equation (\ref{3.3}) has solutions.
\begin{alg}\label{alg3.4}
\begin{itemize}
  \item Construct $T_i$ and $\Lambda_{N_i}$ for each $i\in N.$
  \item Determine $\xi_i\in\R^{k_{N_i}+ k_{-i}}$ as
  \begin{align}\label{3.8}
  \xi_i:=(T_i^\mathrm{T}T_i)^{-1}T_i^\mathrm{T}(V^\phi)^\mathrm{T}.
  \end{align}
  \item Define $\xi_i^1$ as the sub-vector of the first $k_{N_i}$ elements of $\xi_i$. Using (\ref{3.4}), the general form of $V_i, i\in N$ are   calculated as follows
  \begin{align*}
  V_i=(\Lambda_{N_i}\zeta_i+\xi_i^1)^\mathrm{T},
  \end{align*}
  where $\zeta_i$ is an arbitrary vector in $\R^{k_{U(i)}}$.
\end{itemize}
\end{alg}

Using Algorithm \ref{alg3.4}, the local information based utility function of player $i$ is
\begin{align*}
\begin{array}{cll}
c_i(x)&=&V_i\ltimes_{j\in N_i}x_j\\
&=&(\Lambda_{N_i}\zeta_i+\xi_i^1)^\mathrm{T}\mathop\ltimes\limits_{j\in N_i}x_j, \forall i\in N.$$
\end{array}
\end{align*}
\begin{rem}
\begin{enumerate}
\item The computation complexity of Algorithm \ref{alg3.4} is   mainly dependent on the caculation of (\ref{3.8}), where the dimension of $T_i$ is $(k_{N_i}+k_{-i})\times k.$ It shows that the smaller the number of the neighbors is, the lower the computational complexity is.
    Further investigation  for reducing the  computation complexity of Algorithm \ref{alg3.4} is necessary.
\item The method proposed in this paper can also be applied to design state-based potential game \cite{jrm12} with local information based utilities.
\end{enumerate}
\end{rem}
\subsection{Selecting Proper Learning Rule}
After designing the utility function, another thing to be considered is  selecting proper learning algorithm in potential games.  The learning algorithm should ensure that  players' behavior converges to a stategy profile  that maximize the objective function using local information. The logit learning, which is shown as follows, satisfies the above demands:
\begin{itemize}
  \item At each period $t$, player $i\in N$ is chosen with  probability $1/n$ and allowed to update its strategy;
  \item At time $t$, the updating player $i$  selects a strategy $x_i\in S_i$ according to the following probability
$$
{\textbf {Pr}}(x_i(t)=x_i)=\dfrac{{\exp}\{\b{c_i(x_i,x^{-i}(t-1))}\}}{\sum\limits_{y_i\in S_i}{\exp}\{\b{c_i(y_i,x^{-i}(t-1))}\}}
$$
where  $\b\geq0$ is exploration parameter.
\item All other players  repeat their previous actions, i.e., $x^{-i}(t)=x^{-i}(t-1)$, where $x^{-i}(t)\in S^{-i}.$
\end{itemize}
\begin{thm}\cite{ca10} Consider a repeated potential game $G\in{\cal G}_{[n;k_1,\cdots,k_n]}$ with potential function $P(x)$ where all player use  logit learning. The  stationary distribution $\mu^\b\in \D(S)$ is
\begin{align}
\begin{array}{ccl}
\mu^\b(x)=\dfrac{{\exp}\{\b{P(x)}\}}{\sum\limits_{y\in S}{\exp}\{\b{P(y)}\}}
\end{array}
\end{align}
where $\D(S)$ is the set of probability distributions over  $S.$
Moreover, the set of stochastically stable states is equal to the set of maximizers of $P$, where  a state $x\in S$ is called stochastically stable if $\lim_{\b\rightarrow\infty}\mu^\b(x)>0$.
\end{thm}
\begin{rem}
  Since for local information-based potential game we have $c_i(x)=c_i(x_{U(i)})$. Using the above logit learning rule, the local information-based potential game will maximize the potential function $P(x)$ with arbitrarily high probability when $\b$ is sufficiently large. For general potential game if only the local information is allowed to use, the result is not assured. Here general potential game is potential game whose utility functions are not local information based.
\end{rem}
\section{Application: Consensus Problem}

In this section we deal with  consensus problem using the above results. 

Consider a multi-agent system with a system level objective function $\phi:S\rightarrow \R$. The goal of the multi-agent system is to maximum the  objective function. Due to its mobility limitations, agent  can only select strategies from a restricted strategy set. For example, a  robot  in 2-D environment can only move to a position within a radius of its current location. Agent can communicate with its neighbors, which means that at each period the agent can only observe its neighbors' strategies.

A question is: If the player can only select strategies from a restricted strategy set, can the local information-based potential game maximize the objective function  when $\b$ is sufficiently large using logit learning. Unfortunately, the answer is no. But  binary restrictive logit learning, which was introduced in \cite{jm09}, can accomplish the aim. Denote $R_i(x_i(t-1))\subset S_i$ by the  set of strategies available to player $i$ at time $t$. The  binary restrictive logit learning is as follows:
\begin{itemize}
  \item At each period $t$, player $i\in N$ is chosen with  probability $1/n$ and allowed to update its strategy;
  \item At time $t$, the updating player $i$  selects   a trial strategy $x_i^t$ from $R_i(x_i(t-1))$ according to the following probability
  $${\textbf {Pr}}(x_i^t=x_i)=1/w_i, {\text{for~}} x_i\in R_i(x_i(t-1))\setminus x_i(t-1),$$
  $${\textbf {Pr}}(x_i^t=x_i(t-1))=1-((|R_i(x_i(t-1))|-1)/w_i),$$
  where $w_i=\max_{x_i\in S_i}{|R_i(x_i)|}$.
  \item After selecting a trial strategy $x_i^t$, player $i$ chooses its strategy $x_i(t)$ as follows:
\begin{align}\label{4.1}
\begin{array}{ccl}
{\textbf {Pr}}(x_i(t)=x_i^t)=\dfrac{{\exp}\{\b{c_i(x_i^t,x^{-i}(t-1))}\}}{M},\\
{\textbf {Pr}}(x_i(t)=x_i(t-1))=\dfrac{{\exp}\{\b{c_i(x(t-1))}\}}{M},
\end{array}
\end{align}
where $\b\geq0$ is exploration parameter, and
$$M={\exp}\{\b{c_i(x_i^t,x^{-i}(t-1))}\}+{\exp}\{\b{c_i(x(t-1))}\}.$$

\item All other players  repeat their previous actions, i.e., $x^{-i}(t)=x^{-i}(t-1)$, where $x^{-i}(t)\in S^{-i}.$
\end{itemize}

If the restricted strategy sets satisfy the following conditions：
\begin{enumerate}
  \item Reversibility: $$x_i\in R_i(y_i)\Leftrightarrow y_i\in R_i(x_i), \forall x_i,y_i\in S_i, \forall i\in N.$$
  \item Feasibility: For any two strategies $x^1_i,x_i^m\in S_i,$ there exists a series of strategies $x^1_i\rightarrow x^2_i\rightarrow \cdots \rightarrow x_i^m$ such that $x_i^j\in R_i(x_i^{j-1}), j=2,\cdots,m.$
\end{enumerate}
Then using the above binary restrictive logit learning potential game with restricted strategy sets will maximize its potential function $P(x)$ with arbitrarily high probability when $\b$ is sufficiently large \cite{jm09}.
\begin{exa}\label{e4.1}
Consider a multi-agent system with three agents $N=\{1, 2, 3\}$. The system level goal is that all agents gather at point $(3,3)$  in 2-D  environment, such as a room. Each agent has a strategy set $S_i=\{1,2,3\}\times\{1,2,3\}$. Due to its mobility limitations, every agent can only move to a position within a radius $1$ of its current location, and there is an obstacle in $(2,2)$, as shown in Fig.~\ref{Fig4.1}. The communication graph is time-invariant, which is shown in red line. Agent $1$ can only communicate with agent $2$. Agent $2$ can communicate with agent $1$ and $2$. Agent $3$ can only communicate with agent $2$.

\begin{figure}[!hbt]
    \centering
    \includegraphics[height = 5.7cm, width = 6.2cm]{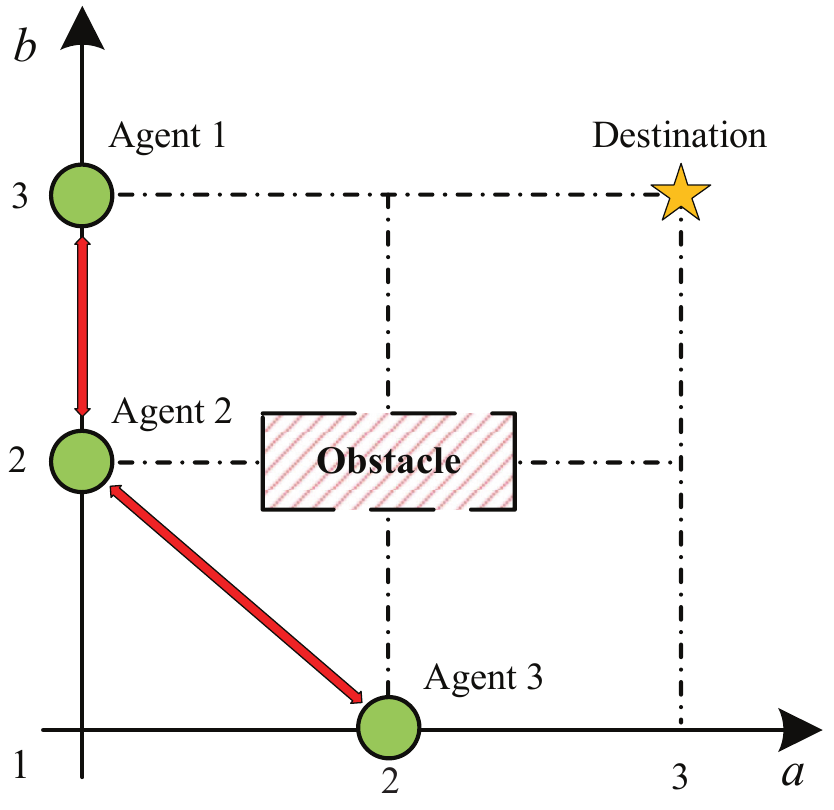}
    \caption{Initial Configuration and Network Graph of Example \ref{e4.1}}
    \label{Fig4.1}  
\end{figure}

 \emph{a). Vector expression of strategies}

Identify
$$\d^i_3\d^j_3\sim(i,j),~~\forall i,j=1,2,3.$$

The system level objective function can be described as
\begin{align}\label{4.34}
\phi(x):=|\{i|x_i(t)=(3,3)\}|,
\end{align}
where $x_i(t)\in S_i, i\in N$. Using the vector expression of strategies, the structure vector of $\phi$ can be calculated as
\begin{align*}
\begin{array}{cll}
V^\phi=\\
\sum\limits_{i=1}^8\sum\limits_{j=1}^8(\d_9^9\otimes\d_9^i\otimes\d_9^j+\d_9^i\otimes\d_9^9\otimes\d_9^j+\d_9^i\otimes\d_9^j\otimes\d_9^9)^\mathrm{T}\\
+\sum\limits_{i=1}^8(\d_9^i\otimes\d_9^9\otimes\d_9^9+\d_9^9\otimes\d_9^i\otimes\d_9^9+\d_9^9\otimes\d_9^9\otimes\d_9^i)^\mathrm{T}~~\\
+(\d_9^9\otimes\d_9^9\otimes\d_9^9)^\mathrm{T}.~~~~~~~~~~~~~~~~~~~~~~~~~~~~~~~~~~~~~~~~~~
\end{array}
\end{align*}

 \emph{b). Local information based utility design}

Set
$$\Gamma=
\left[
  \begin{array}{ccc}
    T_1 & {\bf 0}_{729\times 810}   & {\bf 0}_{729\times 162} \\
    {\bf 0}_{729\times 162}   & T_2 & {\bf 0}_{729\times 162} \\
    {\bf 0}_{729\times 162}   & {\bf 0}_{729\times 810}   & T_2 \\
  \end{array}
\right],$$
where $T_1=[I_{81}\otimes{\bf 1}_9,{\bf 1}_9\otimes I_{81}]$, $T_2=[I_{729},I_{9}\otimes{\bf 1}_9\otimes I_{9}],$ and $T_3=[{\bf 1}_9\otimes I_{81},I_{81}\otimes{\bf 1}_9].$

We can verify that the following equation has  solutions
$$\Gamma\cdot\xi={\bf 1}_3\otimes (V^{\phi})^\mathrm{T}.$$
According to Theorem \ref{th3.1}, the local information based utility function is designable. The designed local information based utility function has the following form:
$$c_i(x)=V_i\ltimes_{j\in N_i}x_j,~i=1,2,3,$$
where $N_1=\{1,2\},~N_2=\{1,2,3\}$, and $N_3=\{2,3\}.$

Using (\ref{3.4}), the general form of $V_i, i=1,2,3$ are   calculated as follows
\begin{align*}
\begin{array}{cll}
V_1^\mathrm{T}&=&({\bf 1}_9\otimes I_9)\cdot\zeta_1+[0.07,0.02,0.02,-0.04,\\
&~&4.8,2.8,\cdots,1,1,2]\in \R^{81};\\
V_2^\mathrm{T}&=&(I_9\otimes{\bf 1}_9\otimes I_9)\cdot\zeta_2+[1.76,-1,-2.2,-1,\\
&~&2.3,3.1,\cdots,1,1,0,0,0,0,0]\in \R^{729};\\
V_3^\mathrm{T}&=&(I_9\otimes{\bf 1}_9)\cdot\zeta_3+[1,1,1,1,1,1,1,1,1,\\
&~&-2.6,-2.6,\cdots,1,1,1,1,1,2]\in \R^{81},
\end{array}
\end{align*}
where $\zeta_1\in\R^{9}, \zeta_2\in\R^{81},$ and $\zeta_3\in\R^{9}$ are arbitrary vectors.

 \emph{c). Simulation results}

The initial configuration of all agents is shown in Fig.~\ref{Fig4.1}. Let $\zeta_1={\bf 1}_{9},\zeta_2={\bf 1}_{81}$, and  $\zeta_3={\bf 1}_{9}.$ Using  binary restrictive logit learning with parameter $\b=0.02t$, we have the following simulation results. As $\b \rightarrow\infty~(t\rightarrow\infty)$, local information based potential game converges to an equilibrium point which maximizes the potential function. Then all agents agree to stay at destination $(3, 3)$ forever, which are shown in Fig. \ref{Fig4.2}-Fig.\ref{Fig4.4}. Here $x_i(t)=(a_i(t),b_i(t)),~i=1,2,3.$
\begin{figure}[!hbt]
    \centering
    \includegraphics[height = 3.7cm, width = 8.6cm]{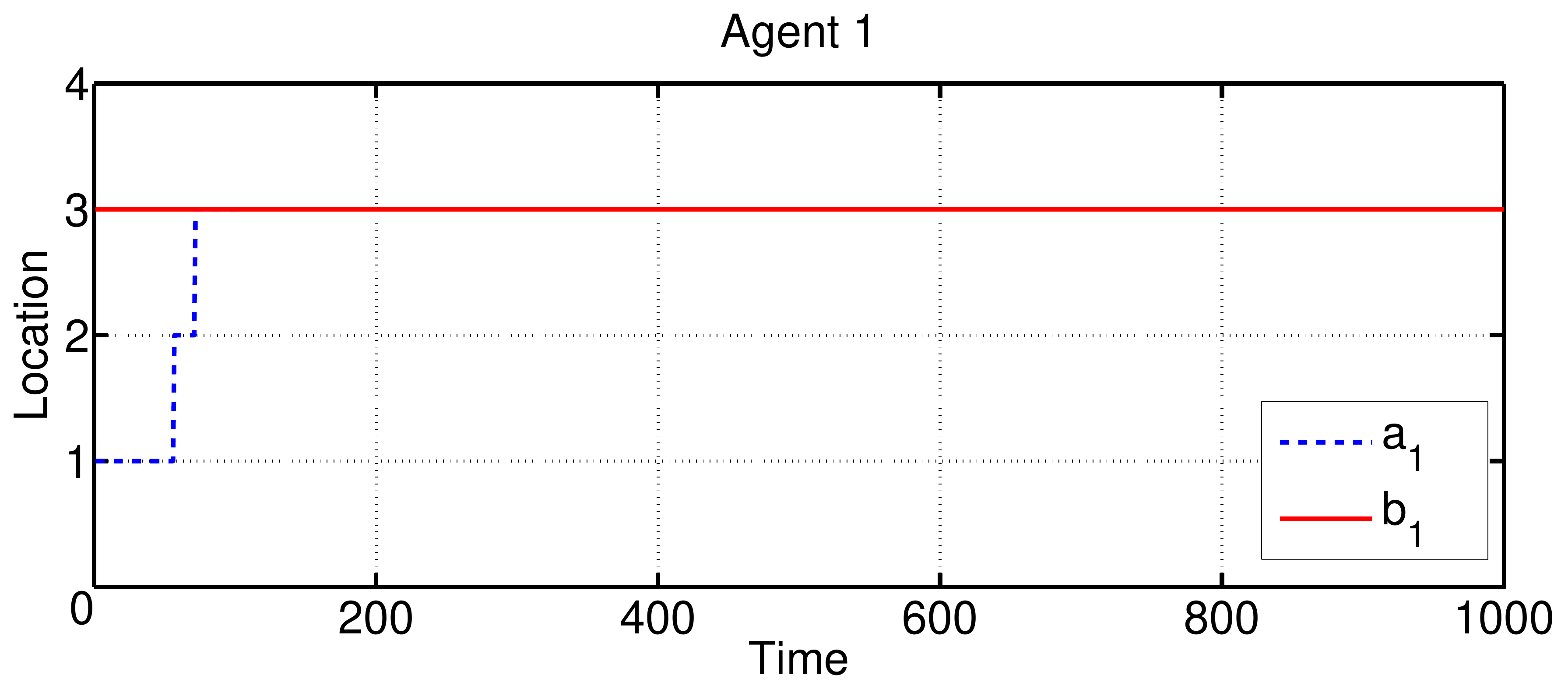}
    \caption{Evolution of Agent 1's Position of Example \ref{e4.1}}
    \label{Fig4.2}  
\end{figure}

\begin{figure}[!hbt]
    \centering
    \includegraphics[height = 3.7cm, width = 8.6cm]{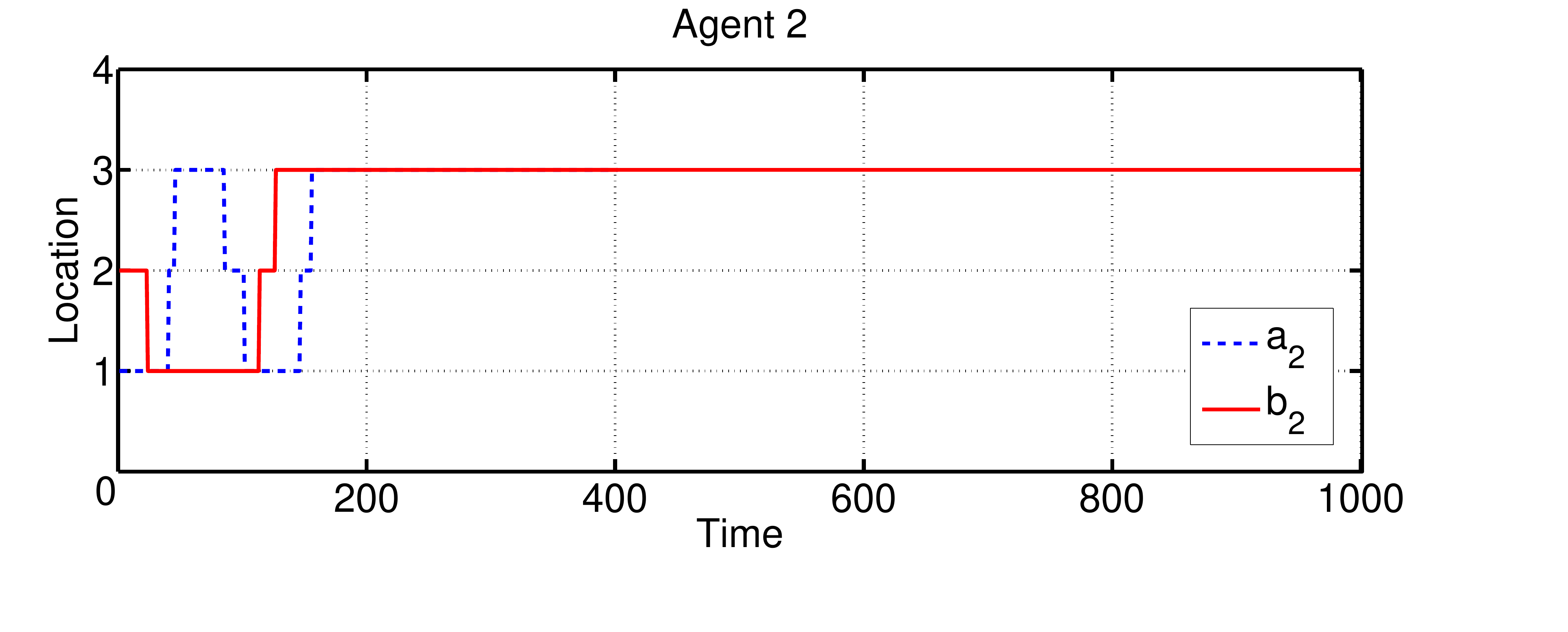}
    \caption{Evolution of Agent 2's Position of Example \ref{e4.1}}
    \label{Fig4.3}  
\end{figure}

\begin{figure}[!hbt]
    \centering
    \includegraphics[height = 3.7cm, width = 8.6cm]{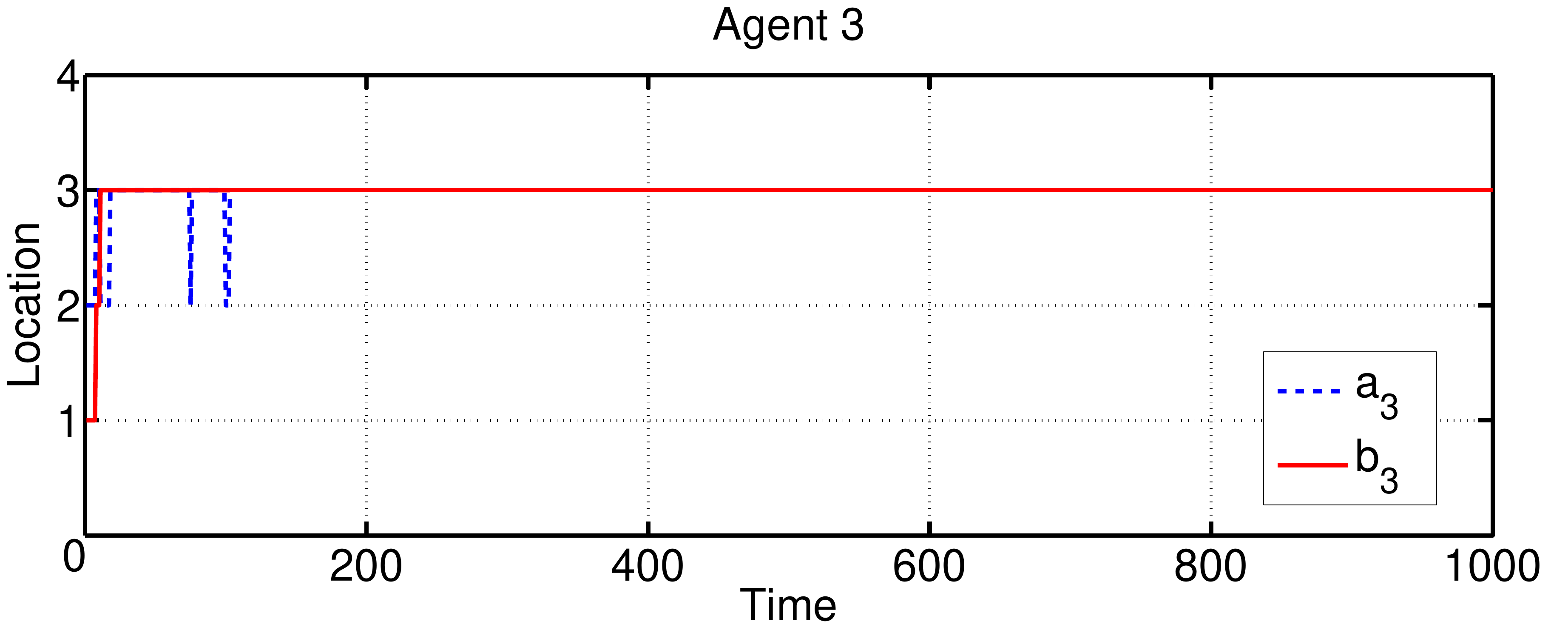}
    \caption{Evolution of Agent 3's Position of Example \ref{e4.1}}
    \label{Fig4.4}  
\end{figure}

\end{exa}


\section{Conclusion}
This paper investigates the design of potential game using local information.
We firstly present  a necessary and sufficient condition for the existence of local information based utility functions, which  can be verified by checking whether a series of linear equations have a solution. Local information based utility functions can be designed using the solutions when the linear equations have solutions. Then a consensus problem is  used to demonstrate the effectiveness of the proposed design  approach.

Open and interesting questions for  further investigations include: If the utility functions cannot be designed as a potential game using local information, can we  design a near-potential game \cite{oc13}? To maximize the system level objective function, how ``near" the designed game should be for a given learning rule?

\end{document}